\documentclass[preprint,12pt]{elsarticle}
\usepackage{amsmath,latexsym,amssymb, enumerate, amsthm, epsfig }
\usepackage{amsmath,amsfonts,amsthm,amscd,amsxtra,enumerate,bm}
\usepackage[paperheight=297mm, paperwidth=210mm, left=30mm, right=20mm]{geometry}
\usepackage{hyperref}
\usepackage[hyperref]{xcolor}
\hypersetup{colorlinks=true,}

\begin{document}

\begin{frontmatter}

\title{Transport equations of nonlinear geometric optics in stratified media exhibiting mixed nonlinearity}


\author{Harsh V. Mahara and V. D. Sharma}

 \address{Department of Mathematics, Indian Institute of Technology Bombay,
	Powai, Mumbai-400076}

\begin{abstract}
In this paper, we concerned with the propagation of sound waves through stratified media. Transport equation of nonlinear geometric optics in media with mixed nonlinearity, in the case of spatially varying density and entropy fields, is derived; this equation contains a cubic nonlinear term in addition to the quadratic nonlinearity. We consider an atmosphere where thermodynamic quantities depend on the height only. Effects of real gas parameters, density, and entropy attenuation parameters on the breakdown of the solution are investigated numerically.

\end{abstract}


\end{frontmatter}


\section{Introduction}

The efficacy of the geometric acoustics method has lead  many researchers to extend the study and develop the  underlying ideas to  wave propagation problems (see, for example Cramer and Sen \cite{MR1166185}, Hunter and Keller\cite{MR716196}, Kluwick and Cox\cite{MR1600236}, Majda and Rosales \cite{MR760229}); the technique requires the introduction of fast and slow variables and the phase functions. Krylov and Bogolibuov \cite{JF85} have employed these methods in the context of ODEs  as a substitute for the method  employed  by Poincar\'e and Lindstedt in the problems of celestial mechanics. The exact scaling of the fast variable with respect to the slow variable may vary depending on the problem under consideration. 

In certain systems, having singular thermodynamic behavior, where the fundamental derivative is small, nonlinear distortions are observed over a time scale longer by an order of magnitude; hence, it is necessary to use a higher order, i.e., $ O(\epsilon^{-2}) $  fast variables. Kluwick and Cox \cite{MR1600236} have used this methodology to the equations of gasdynamics and found that the evolution equation governing the asymptotic behavior contains a quadratic nonlinear term besides a cubic nonlinearity. 

In this chapter, transport equation of nonlinear geometric acoustics in media with mixed nonlinearities following Kluwick and Cox \cite{MR1600236} is derived in the case of stratified fluid and gravitational source terms. The characteristic feature is that cubic nonlinearities arise in addition to delayed nonlinear distortions. We provide singular correction terms to the transport equations that result when the quasilinear equations are manipulated through multiplication by polynomial nonlinearities.
Effectively, any problem of the acoustics occurs in the presence of the gravitational field, and as a result, the unperturbed state ceases to be uniform. For problems considering the propagation of the sound waves over a large distance, such as in the ocean or in the atmosphere, these effects may be significantly important and generate rarefaction and amplification of sound waves.

Considering a quiet and steady atmosphere where thermodynamic quantities depend on the height only and varying according to the exponential laws, a simplified form of the evolution equation is obtained. It is shown that the solution of the Cauchy problem for the evolution equation  exhibits a breakdown of the continuous solution on the expansive phase of the wave profile, which is monotonic increasing, in the sense that the Jacobian of the transformation vanishes after a finite time.  This behaviour is due to the presence of the cubic nonlinearity term in the flux function and is quite different from   the quadratic nonlinearity case where the solution is always continuous.

The work is organized as follows: Basic equations and formulation of the problem are given in Section \ref{short}.  In Section \ref{evo0}, using the ideas of Kluwick and Cox \cite{MR1600236}, a detailed derivation of transport equations is given. The expressions for quadratic and cubic nonlinear terms are obtained in the Section \ref{real}. A brief discussion of the atmospheric model is considered in Section \ref{am}. In Section \ref{real1}, the effects of real gas parameters and atmospheric parameters are seen on the breakdown of the solution. Finally, we conclude this chapter with a discussion of our results in Section \ref{conv}.

\section{Basic equations and short wave limit}\label{short}
Equations describing the propagation of acoustics waves through a stratified fluids may be expressed in the form

\begin{equation}\label{l1.1} 
\left.\begin{aligned}
\frac{\partial \bf u}{\partial t} + {(\bf u.\nabla)u} + \rho ^{-1} \nabla p  = -{\bf G},\\
\frac{\partial \rho}{\partial t} + {(\bf u.\nabla)}\rho + \rho (\nabla. {\bf u}) =0, \\
\frac{\partial s}{\partial t} + {(\bf u.\nabla)}s = 0, \\
\end{aligned}\right.
\end{equation}
where $ \rho $ is the density of the fluid, $ p=p(\rho,s) $ the pressure, $ s $ the entropy, $ t $ the time, $ \nabla $ the gradient operator with respect to the space coordinates $ (x_{1},x_{2},x_{3}) $, and $ {\bf G} $ the forcing function, which balances the initial conditions; a comma followed by the letter $ t $ denotes partial differentiation with respect to time, $ t$. The governing system \eqref{l1.1} can be rewritten into the form 
\begin{equation}\label{l1.2}
\frac{\partial\mathbf{U}}{\partial t}  + A^{k}{\mathbf {(U)}}\frac{\partial\mathbf{U}}{\partial x_{k}}  + {\mathbf{F}}(\mathbf{U})=0,\qquad k=1,2,3,
\end{equation}
representing a quasilinear hyperbolic system of equations with source terms which are attributed to the influences of temperature and gravity. Here $ \mathbf{U} $ and $ \mathbf{F} $ are column vectors defined as $ \mathbf{U}=(u_{1},u_{2},u_{3},\rho,s)^{tr} $ and $ \mathbf{F}=(f_{1},f_{2},f_{3},0,0)^{tr} $, respectively;
$ u_{k} $ are the fluid velocity components and $ f_{k} $ the known function of $ \mathbf{U}_{j} $ and $\mathbf{A}^{k} $ are $ 5 \times 5 $ matrices with entries $ A_{mn}^{k}, 1\leq m,n\leq 5, $
defined as
\begin{equation*}\label{l1.31}
\setlength{\jot}{5pt}
\begin{aligned}
& \quad\qquad 
A_{11}^k = A_{22}^k = A_{33}^k =A_{44}^k =A_{55}^k=u_{k},
\\
&\quad\qquad 
A_{12}^k = A_{21}^k = A_{13}^k=A_{31}^k=A_{23}^k=A_{32}^k=A_{45}^k=A_{54}^k=A_{5j}^k=0,
\\
& \quad\qquad 
A_{4j}^k =\rho \delta_{jk}\quad~ A_{i4}^k = \frac{a^2\delta_{ik}}{\rho},\quad~ A_{i5}^k=\frac{\partial P}{\partial s}\frac{\delta_{ik}}{\rho},
\\
\end{aligned}
\end{equation*}
where $ 1 \leq i,\,j,\,k\leq 3 \;$, $ a^{2} = (\partial p/\partial \rho )|_{s}$ is the speed of sound, and $ \delta $ is the Kroneker delta.\\
We take the unperturbed solution to be $ \mathbf{U}_{0}=(0,0,0,\rho_{0}(x),s_{0}(x))^{tr} $, where the subscript $ 0 $ is used to characterize the unperturbed fluid in equilibrium.\\
We non-dimensionalize the system \eqref{l1.2} by introducing the dimentionless variables which are starred, defined as \\
$ x^{*}=x/L,\quad t^{*}=t\sqrt{gH}/L,\quad u^{*}=u/\sqrt{gH},\quad \rho^{*}=\rho/\rho_{c},\quad p^{*}=p/\rho_{c}gH,\quad a^{*}=a/\sqrt{gH},\\ s^{*}=s/s_{c},\quad{\mathbf{F}}^{*}={\mathbf{F}}/g,  $ \\
where $ L $ is the length of the disturbed region, $ H $ is the scale height of stratification defined as the typical value of $ \rho|\nabla \rho_{0}|^{-1} ,$ $ g $ is the acceleration due to gravity, and $ \rho_{c} $ and $ s_{c} $ are constants representing convenient reference density and entropy, respectively. When $ L $ is much smaller than the scale height $ H $ of stratification, the wave is characterised as a short wave; such waves are common in the atmosphere and ocean.\\
Using nondimensionalize variables in (\ref{l1.2}) gives,\\
\begin{equation}\label{l1.5}
\frac{\partial\mathbf{U}}{\partial t} + A^{k}{\mathbf {(U)}}\frac{\partial\mathbf{U}}{\partial x_{k}}  +\nu(\epsilon) {\mathbf{F}}(\mathbf{U})=0,\qquad k=1,2,3,
\end{equation}
Where $ \nu(\epsilon)=L/H $ with $ \epsilon \ll 1 ,$ so that (\ref{l1.5}) describes short waves in the interior of a stably stratified fluid. In (\ref{l1.5}) and hereafter, we drop the star on the nondimensional variables. In the short wave limit, we have assumed that the perturbations caused by the waves are of size 
$ O(\epsilon), $ and they depend significantly on the fast characteristic variable $\xi=\phi({\mathbf{x}},t)/\epsilon^{2}, $ where $ \phi $ is the phase function to be determined. Using the method of multiple scales, for this change of variables and replaced partial derivatives $ \displaystyle{\frac{\partial}{\partial X} \rightarrow  \frac{\partial}{\partial X} + \epsilon^{-2}\left(\frac{\partial \phi}{\partial X}\right) \frac{\partial}{\partial \xi} },$ ($ X $ being either $ t $ or $ x_{k} $) the system (\ref{l1.5}) becomes 
\begin{equation}\label{l1.6}
\epsilon^{2}\left(\frac{\partial\mathbf{U}}{\partial t} + A^{k}{\mathbf {(U)}}\frac{\partial\mathbf{U}}{\partial x_{k}}  +\nu(\epsilon) {\mathbf{F}}(\mathbf{U})\right) +
\left(\frac{\partial\phi}{\partial t} I + A^{k}{\mathbf {(U)}}\frac{\partial\phi}{\partial x_{k}}\right)\frac{\partial\mathbf{U}}{\partial \xi}=0,
\end{equation}
where $ I $ is the $ 5 \times 5 $ unit matrix. 

\section{Evolution equation}\label{evo0}

We seek for small amplitude high frequency wave solutions of (\ref{l1.2})  with an asymptotic approximation as $ \epsilon \to 0 $ of the form,
\begin{equation}\label{l1.7}
{\mathbf{U}} = {\mathbf{U}_{0}(\mathbf{x})}+\epsilon{\mathbf{U}^{(1)}(\xi,\mathbf{x},t)}+
\epsilon^{2}{\mathbf{U}^{(2)}(\xi,\mathbf{x},t)}+
\epsilon^{3}{\mathbf{U}^{(3)}(\xi,\mathbf{x},t)}+O(\epsilon^{4}),
\end{equation}
which corresponds to the propagation of waves of small amplitude with a disturbed region of small extent. Here $ \epsilon $ is a small parameter measuring the wave amplitude, $ \mathbf{U}_{0}=[0,0,0,\rho_{0}(\mathbf{x}),s_{0}(\mathbf{x})]^{tr} $, the known background state, is a solution of 
\begin{equation}\label{l1.80}
{\mathbf{U}_{0,t}}  + A_{0}^{k}{\mathbf {(U)}}{\mathbf{U}}_{0,x_{k}} =0,
\end{equation}
(recall that $ \nu=0 $ when $ \epsilon=0 $) and $ \mathbf{U}^{(1)}, $$ \mathbf{U}^{(2)} $ and $ \mathbf{U}^{(3)}$ are first, second and third order perturbations of the undisturbed state $\mathbf{U}^{(0)},  $ respectively. We now expand $ A^{k}({\mathbf{U}}) $ and 
$ F({\mathbf{U}}) $ as a power series in $\epsilon $ about $ \mathbf{U} =\mathbf{U}_{0} $ and use them in (\ref{l1.6}) together with (\ref{l1.7})  to obtain the resulting asymptotic expansion in the form
\begin{equation}\label{l1.8}
\epsilon Z_{1}+\epsilon^{2} Z_{2}+\epsilon^{3} Z_{3}+\ldots=0;
\end{equation}
where $ Z_{1}=0 $ is the linearized system associated with (\ref{l1.6}) that admits five families of characteristic surfaces two of which represent waves propagating with speeds $ \pm a_{0} $ through the background state $ \mathbf{U}=\mathbf{U}_{0} $ and the remaining three form a set of coincident characteristics representing entropy waves or particle paths; here we shall be concerned with the propagation of a right runnig acoustic wave $ \phi(x,t)=\text{constatnt}, $ propagating with speed 
$a_{0}|\nabla \phi|,  $ and so the linear rays of (\ref{l1.5}) are characteristic curves (bicharacteristics of (\ref{l1.2}))of 
\begin{equation}\label{l1.9}
\frac{\partial\phi}{\partial t}+a_{0}|\nabla\phi|=0,
\end{equation}
We denote the left and right eigenvectors of $ \sum_{k}(\partial_{k}\phi)A^{k} $ associated with the eigenvalue $ a_{0}|\nabla\phi| $ by $ \mathbf{l} $ and $ \mathbf{r} $, respectively; these are given by\\
\begin{equation}\label{l1.10}
\mathbf{l}=\left(n_{1},n_{2},n_{3},\frac{a_{0}}{\rho_{0}},
(\rho_{0}a_{0})^{-1}\left(\frac{\partial p}{\partial s}\right)_0 \right) \qquad
\mathbf{r}=\left(n_{1},n_{2},n_{3},\frac{\rho_{0}}{a_{0}},
0\right), 
\end{equation} 
where $ \mathbf{n}=\nabla \phi / |\nabla \phi|. $ Furthermore, it follows from $ Z_{1}=0 $ that $ \mathbf{U}^{(1)}=\sigma (\xi ,\mathbf{x},t)\,\mathbf{r} ,$ where $ \sigma  $ is the scaler amplitude function will be determine at the next order. In the present context we refer to the work of Kluwick and Cox \cite{MR1600236} and Cramer and Sen \cite{MR1166185}, who treat the wave propagation problem when the nonlinear effects are noticeable over times of order $ O(\epsilon^{-2}) $ rather than $ O(\epsilon^{-1}); $  their main result focus on the case when the quadratic nonlinearity parameter $ \Gamma $ defined as
\begin{equation}\label{l1.11}
\Gamma= \frac{\partial\phi}{\partial x_{k}}~(\mathbf{l.r})^{-1}~
\mathbf{l}~[\mathbf{r}.(\nabla_{\mathbf{U}}A^{k})_{0}]~\mathbf{r},
\end{equation}
is of order $ O(\epsilon), $ where $ \nabla_{\mathbf{U}} $ is the gradient operator with respect to vector $ \mathbf{U} $ and 
\begin{equation}\label{l1.12}
\mathbf{r.}(\nabla_{\mathbf{U}}A^{k})_{0} = r_{m}\frac{\partial A^{k}}{\partial U_{m}}\Bigr|_{\mathbf{U=U_{0}}}.
\end{equation}
Computation of $ \Gamma $ using the definition of $ A^{k}_{0} $ yield

\begin{equation}\label{l1.13}
\Gamma=\left(1+\frac{\rho_{0}}{a_{0}}\frac{\partial a}{\partial \rho}\Bigr|_{0}\right)|\nabla\phi|,
\end{equation}
which is of order $ O(\epsilon). $ In order to account for a small but nonzero $ \Gamma ,$ it would be convenient to write (\ref{l1.8}) as
\begin{equation}\label{l1.14}
\epsilon Z_{1}+\epsilon^{2} (Z_{2}-\mu)+\epsilon^{3} (Z_{3}+\bar{\mu})+\ldots=0,
\end{equation}
where $ \mu=\Gamma\,\sigma\, \frac{\partial \sigma}{\partial \xi}\,\mathbf{r}$ and  $ \bar{\mu}=\Gamma\,\mathbf{r}\,\sigma\, \frac{\partial \sigma}{\partial \xi}/\epsilon $, noting that $ \Gamma/\epsilon= O(1). $ Thus to the second order, we obtain $ (Z_{2}-\mu) = 0 $ or equivalently,
\begin{equation}\label{l1.15}
\left(\frac{\partial\phi}{\partial t} I + A^{k}{\mathbf {(U)}}\frac{\partial\phi}{\partial x_{k}}\right)\frac{\mathbf{\partial U}}{\partial\xi}^{(2)} =
-\left\lbrace \frac{\partial \phi}
{\partial x_{k} }[\mathbf{r}.(\nabla_{\mathbf{U}}A^{k})_{0}]\mathbf{r}-
\Gamma\mathbf{r}\right\rbrace \sigma\frac{\partial \sigma}{\partial \xi},
\end{equation}
Thus the solvability condition for $ \mathbf{U}^{(2)}, $ which requires that the right hand side of (\ref{l1.15}) be orthogonal to $ \mathbf{l} ,$ is satisfied automatically in view of the definition of $ \Gamma .$ In order to have a solution which exhibit the character of a progressive wave in which both nonlinear and source terms are present, we must get, to the third order, a new equation involving the source terms; for this one must choose $ \nu(\epsilon) $ of the order $ O(\epsilon) $ to keep the source term in 
$ Z_{3}+\bar{\mu}=0. $ Thus, choosing $ \nu(\epsilon) = \epsilon,$ we obtain to the third order;

\begin{equation}\label{l1.16}
\begin{split}
\left(\frac{\partial\phi}{\partial t} I + A^{k}{\mathbf {(U)}}\frac{\partial\phi}{\partial x_{k}}\right)
&
\frac{\mathbf{\partial U}}{\partial\xi}^{(3)}=
\frac{\partial (\sigma \mathbf{r}) }{\partial t}-
A_{0}^{k}\frac{\partial (\sigma \mathbf{r} )}{\partial x_{k}} -[\mathbf{r}.(\nabla_{\mathbf{U}}A^{k})_{0}]
\frac{\mathbf{\partial U}}{\partial x_{k}}^{(0)}\sigma-\\
&
-\left(\frac{\partial\phi}{\partial x_{k}}
[\mathbf{U}^{(2)}.(\nabla_{\mathbf{U}}A^{k})_{0}]\mathbf{r}\sigma_{,\xi}
+\frac{1}{2}\frac{\partial  \phi}{\partial x_{k}}
[\mathbf{r}\mathbf{r}:(\nabla_{\mathbf{U}}\nabla_{\mathbf{U}}A^{k})_{0}]
\mathbf{r}\sigma^{2}\right)\frac{\partial \sigma}{\partial \xi}\\
&
\widehat{\Gamma}\sigma\frac{\partial \sigma }{\partial t}\mathbf{r}
-\frac{\partial\phi}{\partial x_{k}}
[\mathbf{r}.(\nabla_{\mathbf{U}}A^{k})_{0}]
\frac{\mathbf{\partial U}}{\partial\xi}^{(3)}\sigma
-\mathbf{F}_{0}-
[\mathbf{r}.(\nabla_{\mathbf{U}}\mathbf{F})_{0}]\sigma,
\end{split}
\end{equation}\\
where  $ \quad\mathbf{r}\mathbf{r}:(\nabla_{\mathbf{U}}\nabla_{\mathbf{U}}A^{k})_{0}=\sum_{m,n=1}^{5}r_{m}\,r_{n} \frac{\partial^{2}A^{k}}{\partial U_{m}\,\partial U_{n}}|_{\mathbf{U=U_{0}}},\quad
\text{and}~\, \widehat{\Gamma} = \Gamma/\epsilon = O(1). $

Solvability condition for $ \mathbf{U^{(3)}} $ yields the desired evolution equation for $ \sigma, $ namely 
\begin{equation}\label{l1.17}
\frac{\partial\sigma}{\partial\tau}+
\left(\hat{\Gamma}+\frac{E}{2}\sigma\right)\sigma\frac{\partial\sigma}{\partial\xi}+
\left(\mathbf{M}.\mathbf{U}^{(2)}\right)\frac{\partial \sigma}{\partial \xi}+
\left(\mathbf{N}.\frac{\partial\mathbf{U}^{(2)}}{\partial \xi}\right) \sigma+ 
\chi \sigma + \frac{\mathbf{l}.\mathbf{F}_{0}}{2}+
\frac{1}{2}\frac{\Gamma(\mathbf{l}.\mathbf{U}^{(2)}\sigma)}{\partial \xi}=0,
\end{equation}
where $ \partial /\partial \tau =\partial/\partial t + a_{0}\,\mathbf{n}.\nabla $ is  the ray derivative with $ \nabla $ as the gradient operator with respect to the space variable $ x_{k} $, and $ E ,$ 
$ \mathbf{b}, $ $ \mathbf{c} $ and $ \chi $ are given by
\begin{equation}\label{l1.18}
E=\frac{1}{2} \frac{\partial\phi}{\partial x_{k}}
\mathbf{l}\,.[\mathbf{r}\,\mathbf{r}:(\nabla_{\mathbf{U}}\nabla_{\mathbf{U}}A^{k})_{0}]\,\mathbf{r},\qquad\qquad\qquad\qquad\qquad\qquad\qquad\qquad
\end{equation}
\begin{equation}\label{l1.19}
\mathbf{M}=\frac{1}{2}\left\lbrace\frac{\partial \phi}{\partial x_{k}}
\mathbf{l}\,[(\nabla_{\mathbf{U}}A^{k})_{0}]\mathbf{r}-\Gamma\,\mathbf{l}\right\rbrace,\qquad\qquad\qquad\qquad\qquad\qquad\qquad\quad
\end{equation}
\begin{equation}\label{l1.20}
\mathbf{N}=\frac{1}{2}\left\lbrace\frac{\partial \phi}{\partial x_{k}}
\mathbf{l}\,[(\mathbf{r}.\nabla_{\mathbf{U}}A^{k})_{0}]-\Gamma\,\mathbf{l}\right\rbrace,\qquad\qquad\qquad\qquad\qquad\qquad\qquad\quad
\end{equation}
\begin{equation}\label{l1.21}
\mathbf{\chi}=\frac{1}{2}\left\lbrace
\mathbf{l}\;A_{0}^{k}.\frac{\partial\mathbf{r}}{\partial x_{k}} +
\mathbf{l}\;[\mathbf{r}.(\nabla_{\mathbf{U}}A^{k})_{0}]\;
\frac{\mathbf{\partial U}}{\partial x_{k}}^{(0)}+
\mathbf{l}\;[\mathbf{r}.(\nabla_{\mathbf{U}}\mathbf{F})_{0}]	\right\rbrace.\qquad\qquad\qquad\qquad\qquad
\end{equation}
It may be notice that in (\ref{l1.17}), it will be essential to have a knowledge of the second order approximation, $ \mathbf{U}^{(2)} $ and
$ \mathbf{U}^{(2)}_{,\xi} $ ; this can be achieved by noting that (\ref{l1.15}) can be integrated to yield
\begin{equation}\label{l1.22}
\mathcal{A}.\mathbf{U}^{(2)}=
-\frac{1}{2}\left(
\frac{\partial \phi}{\partial x}\left[
\mathbf{r}.(\nabla_{\mathbf{U}}A^{k})_{0}
\right]
\mathbf{r}-\Gamma\mathbf{r}
\right)\sigma^{2},\qquad\qquad\qquad\qquad\qquad\qquad\qquad\quad
\end{equation}
where $ \mathbf{U}^{(2)}=0 $ when $ \mathbf{U}^{(1)}=0, $ and $ \mathcal{B}=[\frac{\partial \phi}{\partial t}I+A_{0}^{k}\frac{\partial \phi}{\partial x_{k}}]. $ In view of Eqs. $(\ref{l1.13}), $ $ (\ref{l1.19}), $ and $ (\ref{l1.20}),$ we find that $ \mathbf{M} $ and $ \mathbf{N} $ are orthogonal to $ \mathbf{r} $ and hence they lie in the linearly independent row space of $ \mathcal{B} .$ Hence, we can write
\begin{equation}\label{l1.23}
\mathbf{M}=\mathcal{\omega}_{\alpha}~\mathcal{B}_{\alpha},\qquad
\mathbf{N}=\mathcal{\delta}_{\alpha}~\mathcal{B}_{\alpha},\qquad\qquad\qquad\qquad\qquad\qquad\qquad\qquad\qquad\qquad\qquad
\end{equation}
where  $ \mathcal{B}_{\alpha}~(\alpha=1,2,3,4)$ being the linearly independent rows of $ \mathcal{B}$  
then, the terms $ \mathbf{M}.\mathbf{U}^{(2)} $ and $ \mathbf{N}.\mathbf{U}^{(2)}_{,\xi} $ in Eq. (\ref{l1.17}) can be written as
\begin{equation}\label{1.24}
\mathbf{M}.\mathbf{U}^{(2)}=-\frac{\bm{\omega}}{2}\left(
\frac{\partial \phi}{\partial x}\left[\mathbf{r}.(\nabla_{\mathbf{U}}A_{1}^{k})_{0}\right]-\Gamma I_{1}\right)\mathbf{r}\sigma^{2},\qquad\qquad\qquad\qquad\qquad\qquad\qquad\qquad
\end{equation}
\begin{equation}\label{1.25}
\mathbf{N}.\mathbf{U}^{(2)}_{\xi}=-\bm{\delta}\left(
\frac{\partial \phi}{\partial x}\left[\mathbf{r}.(\nabla_{\mathbf{U}}A_{1}^{k})_{0}\right]-\Gamma I_{1}\right)\mathbf{r}\sigma\sigma_{,\xi},\qquad\qquad\qquad\qquad\qquad\qquad\qquad\qquad
\end{equation}
where $ A_{1}^{k} $ and $ I_{1} $ being $ 4\times5 $ matrices obtained from $ A^{k} $ and 
$ I $, respectively by deleting their last row,
$ \bm{\omega }=(\omega_{1},\omega_{2},\omega_{3},\omega_{4}) $ and  
$ \bm{\delta}=(\delta_{1},\delta_{2},\delta_{3},\delta_{4},). $ Thus the evolution equation (\ref{l1.17}) becomes
\begin{equation}\label{l1.26}
\frac{\partial\sigma}{\partial\tau}+
\left(\hat{\Gamma}+\frac{\Lambda}{2}\sigma\right)\sigma\frac{\partial\sigma}{\partial\xi}+
\chi \sigma+
\frac{\mathbf{n}.\mathbf{F}_{0}}{2}=0,
\end{equation}
where
\begin{equation}\label{l1.27}
\Lambda= E-(\bm{\omega}+2\bm{\delta})\left(\frac{\partial \phi}{\partial x}\left[\mathbf{r}.(\nabla_{\mathbf{U}}A_{1}^{k})_{0}\right]-\Gamma I_{1}\right)\mathbf{r}.\qquad\qquad\qquad\qquad\qquad\qquad\qquad
\end{equation}
We now turn to the calculation of coefficients $ \hat{\Lambda}, $ $\Gamma  $ and $ \chi .$
If we define 
$ \nabla_{l}A_{ij}^{k}=\left(\frac{\partial A_{ij}^{k}}
{\partial U_{l}}\right)|_{\mathbf{U}=\mathbf{U}_{0}} $ and 
$ \nabla_{lm}A_{ij}^{k}=\left(\frac{\partial^{2} A_{ij}^{k}}
{\partial U_{l}\partial U_{m}}\right)|_{\mathbf{U}=\mathbf{U}_{0}}, $
we find that
\begin{equation}\label{l1.28}
\nabla_{k} A_{ij}^{k}=1, \qquad 
\nabla_{4} A_{k4}^{k}= -\frac{a_{0}^{2}}{\rho_{0}^{2}}+\frac{2a_{0}}{\rho_{0}}a_{\rho0},\qquad
\nabla_{5} A_{k4}^{k}= \frac{2a_{0}}{\rho_{0}}a_{s0},\qquad\qquad\qquad\qquad
\end{equation}
\begin{equation}\label{l1.29}
\nabla_{5} A_{k5}^{k}=\frac{p_{,ss0}}{\rho_{0}},\qquad 
\nabla_{4} A_{k5}^{k}= -\frac{p_{,s0}}{\rho_{0}^{2}}+\frac{2a_{0}}{\rho_{0}}a_{\rho0},\qquad 
\nabla_{4} A_{4k}^{k}=1; \quad 1\leq i \leq 5,\; k=1,2,3,
\end{equation} 
and hence,
\begin{equation}\label{l1.30}
\Gamma = \left(1+\frac{\rho_{0}}{a_{0}}a_{,\rho0}\right)|\nabla\phi|,\qquad\qquad\qquad\qquad\qquad\qquad\qquad\qquad\qquad\qquad\qquad\qquad
\end{equation}
As $\, \Gamma=O(\epsilon), $  the term containing  $ \Gamma $ and $ \Gamma^{2} $ can be neglected for evaluating  $ \Lambda $ and $ \chi .$ Further, since $ \nabla_{lm}A_{ij}^{k} $ in $ E $ appears in the evolution equation as a combination $ r_{j} ,$ $ r_{l} $ and $ r_{m}, $ it is easily verified that all the terms $ \nabla_{lm}A_{ij}^{k} $ vanish except $ \nabla_{44}A_{k4}^{k} $ which is given by
\begin{equation}\label{1.31}
\nabla_{44}A_{k4}^{k}=\frac{2a_{0}^{2}}{\rho_{0}^{3}}(6+\Omega).\qquad\qquad\qquad\qquad\qquad\qquad\qquad\qquad\qquad\qquad\qquad\qquad
\end{equation}
Hence, $ E=(6+\Omega)|\nabla\phi|/a_{0} $ with 
$ \Omega =\frac{\rho_{0}^{2}}{a_{0}}\frac{\partial\Sigma}{\partial\rho}(\rho_{0},s_{0})=O(1) $ and 
$ \Sigma=\frac{1}{\rho}\frac{\partial(a\rho)}{\partial\rho}.$ From (\ref{l1.19}) and (\ref{l1.20}),the vectors 
$ \mathbf{M} $ and $ \mathbf{N} $ are obtained as follows:
\begin{equation}\label{1.32}
\mathbf{M}=\mathbf{N}=|\nabla\phi|\left(n_{1},n_{2},n_{3},-\frac{a_{0}}{\rho_{0}},a_{,s0}\right).\qquad\qquad\qquad\qquad\qquad\qquad\qquad\qquad\qquad
\end{equation}
Using (\ref{l1.23}), we find that
\begin{equation}\label{l1.33}
\omega_{\alpha}=\delta_{\alpha}=\frac{n_{\alpha}}{a_{0}}\left(\frac{a_{0}\rho_{0}a_{,s0}}{p_{,s0}}\right),\quad \alpha=1,2,3; \quad
\omega_{4}=\delta_{4}=\frac{1}{\rho_{0}}\left(\frac{a_{0}\rho_{0}a_{,s0}}{p_{,s0}}+1\right),\quad
\end{equation}
and hence, 
$ (\bm{\omega}+2\bm{\delta})\left(\phi_{,k}\left[\mathbf{r}.(\nabla_{\mathbf{U}}A_{1}^{k})_{0}\right]-\Gamma I_{1}\right)\mathbf{r}=6|\nabla \phi|/a_{0}; $ taking the vector $ \mathbf{g} $ as the acceleration due to gravity and using the forgoing results, we find that 
\begin{equation}\label{l1.34}
\Lambda=\Omega |\nabla \phi|/a_{0}\qquad\text{and}\qquad \chi= \left(a_{0}\nabla.\mathbf{n}+a_{,s0}s_{0,k}n_{k}\right)/2,
\qquad\qquad\qquad\qquad\quad
\end{equation}
where $ \nabla.\mathbf{n} $ is the mean curvature of the wavefront.\\
It is noticeable that the quadratic nonlinearity coefficient $ \hat{\Gamma} $ in (\ref{l1.26}) is  Lax \cite{Lax}  genuine nonlinearity coefficient, whereas the cubic nonlinearity coefficient   $ \Lambda $, indicates the degree of material nonlinearity; the source term $ \chi\sigma $ in (\ref{l1.26})  corresponds to the changes in wave amplitude  $ \sigma $ attributed to the wave interactions with the changing medium ahead and the wavefront curvature as the wave moves along the rays. 
\section{Real gas parameters}\label{real}
For the van der Waals gas, where the proper volume of the gas molecule is reduced by an amount $ \beta  $ and the gas pressure is reduced with respect to the ideal pressure due to the attractive interaction of the molecule, the resulting equation of the state is \cite{zhao2011admissible}
\begin{equation}\label{l1.35}
(p+\alpha \rho^{2})(1-\beta\rho) =\rho R T
\end{equation}
where $ T$ is the temperature, $ R $ the gas constant, and the constants $ \alpha $ and $ \beta $ depends on a particular gas. The expression for the entropy can be obtained from the first law of thermodynamics as 
$ \rho^{\gamma}exp((s-s_{0})/c_{v})=(p+\alpha\rho^{2})(1-\beta\rho)^{\gamma}, $ where $ c_{v} $ is the specific heat at constant volume and $ \gamma $ is the specific heat ratio. The speed of sound 
$ a^{2}=\left(\frac{\partial p}{\partial \rho}\right)^{(1/2)} $is thus given by
$ a^{2}=\frac{\gamma(p+\alpha\rho^{2})}{\rho(1-\beta\rho)}-2\rho\alpha. $ Then the quadratic nonlinearity parameter $ \Gamma $
turns out to be
\begin{equation}\label{l1.36}
\Gamma =\left(1+\frac{\rho_{0}}{a_{0}}a_{,\rho0}\right)|\nabla\phi|= \frac{(\gamma+1)}{2}\left(\frac{1}{(1-\beta\rho_{0})}-
\frac{2\alpha\rho_{0}(2-\gamma-3\beta\rho_{0})}{a_{0}^{2}(\gamma+1)(1-\beta\rho_{0})}\right)|\nabla\phi|,
\end{equation}

Further if we choose $ \alpha  $ and $ \beta $ such that 
$\left(\frac{1}{(1-\beta\rho_{0})}-\frac{2\alpha\rho_{0}(2-\gamma-3\beta\rho_{0})}
{a_{0}^{2}(\gamma+1)(1-\beta\rho_{0})}\right)|\nabla\phi| = \epsilon,$ then $ \Gamma= O(\epsilon) $ whilest 
$ \hat{\Gamma}=\Gamma/\epsilon=(\gamma+1)/2=O(1). $ \\
Similarly we can get the expression of cubic nonlinearity parameter, neglecting terms of the form $ \Gamma $ and $ \Gamma^{2} $ we get,
\begin{equation*}
\Omega=\frac{\rho_{0}^{2}}{a_{0}}\frac{\partial}{\partial\rho}
\left(\frac{1}{\rho}\frac{\partial(a\rho)}{\partial\rho}\right)(\rho_{0},s_{0})=
-\left(\frac{3(1+\gamma)}{2(1-\beta\rho_{0})}-\frac{3\alpha\beta\rho_{0}^{2}}{(1-\beta\rho_{0})a_{0}^{2}}\right)=O(1)
\end{equation*}
neglecting $ O(\beta^{2}) $ terms the expression for $ \Lambda $ can be found from the relation
\begin{equation}\label{l1.37}
\Lambda=\frac{\Omega |\nabla \phi|}{a_{0}}=
-\left(\frac{3(1+\gamma)}{2a_{0}}(1+\beta\rho_{0})-\frac{3\alpha\beta\rho_{0}^{2}}{a_{0}^{3}}\right)|\nabla \phi|,
\end{equation}
and an expression used in the calculation of $ \chi  $ can also be found 
$ \left(\frac{a_{,s0}}{p_{,s0}}\right)= \gamma a_{0}\rho_{0}(1-\beta\rho_{0})^{-1}.$
\section{Atmospheric model}\label{am}
We now consider the wave propagations in the troposphere region of the earth's atmosphere ( approximately $ 0-12 $ km above the earth's surface) with a quiet and steady atmosphere with thermodynamic quantities depending on height only, i.e., they depend on one coordinate only, say $ x_{3} ,$ and satisfy the exponential laws based on the U.S. Standard Atmosphere ($ 1966 $ supplement) (\cite{NZ73}) ,
\begin{equation}\label{l1.38}
\rho_{0}= \exp(-\theta x_{3})\quad \text{and}\quad a_{0}= \exp(-\omega x_{3}),
\end{equation}
where $ \theta,\omega\geq0 $ are attenuation rates for the density and sound speed, respectively. The dependence of the entropy on $ x_{3} $ can be obtained from the equilibrium condition $\nabla p |_{0}=0 ,$ so that
\begin{equation}\label{l1.39}
\frac{\partial a}{\partial s}\Bigr|_{0} s_{0,k}n_{k}=
\left(\frac{a_{,s0}}{p_{,s0}}\right)(-a_{0}^{2}\rho_{0}^{'})=
\frac{\gamma \theta \exp(-\omega x_{3})}{2(1-\beta \exp(-\theta x_{3}))},
\end{equation}
where prime denotes derivative with respect to $ x_{3} ;$ \eqref{l1.39} will be needed in the computation of $ \chi $ which appears in the transport equation \eqref{l1.26}. If we choose our initial disturbance on the horizontal plane 
$ x_{3}=0, $ the characteristic surface (or wavefront) at any time $ t $ can be obtained by solving the eikonal equation
\eqref{l1.9}, indeed the characteristic rays are the solution of the ODE's
\begin{equation}\label{l1.40}
\dfrac{dt}{d\tau}=1, \quad
\dfrac{d\mathbf{x}}{d\tau}=\dfrac{\nabla \phi}{| \nabla\phi|}a_{0}(x_{3}),\quad
\dfrac{d(\partial\phi/\partial x_{i})}{d\tau}=-|\nabla \phi|a_{0}^{'}(x_{3}),\quad
\end{equation}
let $\, \tau=0, $ at $\, t=0, $ so that (\ref{l1.20}) implies that $ \tau= t $ along the rays. Since 
$ \nabla\phi|_{t=0}=(0,0,1), $ it follows from (\ref{l1.38}) and  (\ref{l1.39}), that at any time $ t $ the vector 
$ \nabla\phi $ and the location of the wave front is given by 
\begin{equation}\label{l1.41}
\nabla\phi=(0,0,1+\omega t)\quad \quad x_{3}=\omega^{-1}ln(1+\omega t),\quad\quad
\end{equation}
which describes an ascending wave as $ x_{3} $ is an increasing function of time. The source term $ \mathbf{F} $ in (\ref{l1.2}) is indeed $ \mathbf{F}=(0,0,g,0,0) ,$ which in view of the dimentionless quantities assume the form 
$ \mathbf{F}=(0,0,1,0,0) $ and thus the coefficient $ \hat{\Gamma },$ $ \Lambda ,$ $ \chi $ and the inhomogeneous term $ \mathbf{n.F_{0}} $ in (\ref{l1.26}) are now explicitely known thus finally the evolution equation for the amplitude 
$ \sigma $ describes the propagation of an acoustic wave in a stratified medium becomes,

\begin{eqnarray}\label{l1.42}
\frac{\partial\sigma}{\partial t} + 
\left(\frac{\gamma +1}{2}\right)(1+\omega t)\sigma \frac{\partial\sigma}{\partial \xi}
-\frac{3}{4}(\gamma+1)(1+\omega t)^{2}\biggl(1+\beta(1+\omega t)^{(-\theta/\omega)}-\qquad\qquad\nonumber\\
2\alpha\beta(1+\omega t)^{(2-2\theta/\omega)} \biggr)\sigma^{2}\frac{\partial\sigma}{\partial \xi}+
\frac{\gamma \theta (1+\omega t)^{(\theta/\omega )-1}\sigma}{4((1+\omega t)^{\theta}-\beta)}+\frac{1}{2}=0,\qquad
\end{eqnarray}

by the method of characteristics, the solution along the characteristics 
\begin{equation}\label{l1.43}
\dfrac{d\xi}{dt}=\left(\frac{\gamma +1}{2}\right)(1+\omega t)\sigma-
\frac{3}{4}(\gamma+1)(1+\omega t)^{2}\left[1+\beta(1+\omega t)^{(-\theta/\omega)}-
2\alpha \beta(1+\omega t)^{2(1-\theta/\omega)} \right]\sigma^{2}
\end{equation} is given by
\begin{equation}\label{l1.44}
\dfrac{d\sigma}{dt}=-\frac{\gamma \theta (1+\omega t)^{(\theta/\omega )-1}\sigma}{4((1+\omega t)^{\theta}-\beta)}-\frac{1}{2},
\qquad\qquad\qquad\qquad\qquad\qquad\qquad
\end{equation}
which is an ODE and can be solved to obtain
\begin{equation}\label{l1.45}
\begin{aligned}
\sigma(\xi,t)=&\left(\sigma_{0}+\frac{2(1-u^{1\frac{-\gamma \theta}{4\omega}})}
{(\gamma \theta +4 \omega)}\right)u^{\frac{-\gamma \theta}{4\omega}}+\\
&\frac{\gamma\beta}{4}\left(\sigma_{0}u^{\frac{-\gamma \theta}{4\omega}}(u^{\frac{-\theta}{\omega}}-1)-
\frac{2(u^{1-\frac{\theta}{\omega}}-u^{\frac{-\theta}{\omega}(1+\frac{\gamma}{4})})}{4\omega+\gamma\theta}+
\frac{2(u^{1-\frac{\theta}{\omega}}-u^{\frac{-\gamma\theta}{4\omega}})}{4\omega+\theta(\gamma -4)}
\right)\qquad\qquad
\end{aligned}
\end{equation} 
where $ \sigma_{0} $ is the initial value of $ \sigma $ and $ u=(1+\omega t) $ using this value of $ \sigma $ in characteristic equation we can get a condition for the Jacobian. If we choose initial data at time $ t=t_{0} $ as
\begin{equation}\label{examplesin}
\sigma(\xi,t_{0})=
\left\{ \begin{array}{l l}
sin(\xi) & \quad\ \xi\in (0,\pi/2),\\
~~0 & \quad \text{~otherwise}
\end{array} \right.
\end{equation}
expression for the Jacobian can be writen as
\begin{align}
\xi_{\eta}=\frac{(u^{6-\theta/\omega(2+\gamma/4)}-1)}{24\omega-\theta(8+\gamma)}
\left(\frac{-24\alpha\beta\cos(\eta) }{(\gamma \theta + 4 \omega)}\right)+
\frac{(u^{5-\theta/\omega(2+\gamma/4)}-1)}{20\omega-\theta(8+\gamma)}
\left(\frac{24\alpha\beta\cos(\eta) }{(\gamma \theta + 4\omega)}+12\alpha\beta\cos(\eta)sin(\eta)\right)\nonumber\\
-\frac{(u^{4-\theta/\omega(2+\gamma/4)}-1)}{16\omega-\theta(8+\gamma)}
\left(\frac{\gamma}{(\gamma (\theta-4) + 4\omega)}-\frac{2(2+\gamma)}
{(\gamma\theta + 4\omega)}\right)(1+\gamma)\beta cos(\eta)\qquad\qquad\qquad\qquad\qquad\qquad\nonumber\\
-\frac{(u^{3-\theta/\omega(1+\gamma/2)}-1)}{6\omega-\theta(2+\gamma)}
\left(\frac{\gamma}{(\gamma (\theta-4) + 4\omega)}-\frac{2(2+\gamma)}
{(\gamma\theta + 4\omega)}\right)(1+\gamma)\beta cos(\eta)\qquad\qquad\qquad\qquad\qquad\qquad\nonumber\\+
\frac{(u^{2-\theta/\omega(1+\gamma/2)}-1)}{4\omega-\theta(2+\gamma)}
\left(\frac{(1+\gamma)\gamma\beta}{8}\right)cos(\eta)+
\frac{(u^{4-\gamma\theta/4\omega}-1)}{16\omega-\gamma\theta}
\left(\frac{(4-\gamma\beta)}{\gamma\theta+4 \omega}\right)3(\gamma +1)cos(\eta)
\qquad\qquad\nonumber\\
+\frac{(u^{3-\gamma\theta/4\omega}-1)}{12\omega-\gamma\theta}
\left(\frac{\gamma\beta}{(\gamma-4)\theta+4 \omega}\right)3(\gamma +1)cos(\eta)
+\frac{(u^{2-\gamma\theta/4\omega}-1)}{8\omega-\gamma\theta}
\left(1-\frac{\gamma \beta}{4}\right)2(\gamma +1)cos(\eta)\quad\nonumber\\
-\frac{(u^{3-\gamma\theta/2\omega}-1)}{6\omega-\gamma\theta}
\left(\frac{(4-\gamma\beta)}{\gamma\theta+4 \omega}+2sin(\eta)\left(1-\frac{\gamma\beta}{2}\right)\right)
\left(\frac{3(1+\gamma)cos(\eta)}{2}\right)+1.\quad\quad\quad\quad\quad\quad\quad\qquad
\end{align}
Now we consider two cases for seeing the effect of parameter $ \theta $ and $ \omega  $ on breaking of solution.\\
$ \bullet  $ If we consider the case when density is constant, i.e., $ \theta= 0 $ and taking the initial condition as above the expression for jacobian reduces into the following form

\begin{align}
t^{6}[(-\omega/2 )\alpha \beta cos(\eta) ]+ 
t^{5}[(6/5)(\omega cos(\eta)-2 )\omega^{3}\alpha \beta ] +
t^{4}[((1+\gamma)(1+\beta)/8)-(3/2)\alpha \beta +\nonumber\\
2sin(\eta)\alpha \beta \omega ]3\omega^{2}cos(\eta) +
t^{3}[(1+\gamma)(1+\beta)(1-\omega sin(\eta))-4\alpha \beta (1-3\omega sin(\eta)) ]\omega cos(\eta) +\nonumber\\
t^{2}[3(1+\gamma)(1+\beta)((1/4)-\omega sin(\eta))-3\alpha \beta ((1/2)-4\omega sin(\eta) +\omega (1+\gamma)/2]cos(\eta)+\nonumber\\
t[6sin(\eta)(\alpha \beta -(1+\gamma)(1+\beta)/2) +(1+\gamma)]cos(\eta)+
2=\xi_{\eta}.\qquad\qquad\qquad\qquad\qquad\qquad
\end{align}
$ \bullet $ Similarly when we consider the case when speed of sound parameter is zero, i.e.,
$ \omega \to 0 $  and considering the expression $ \lim \omega\to 0 (1+\omega t)^{(1-\gamma \theta/4\omega )}= 
e^{-(\gamma \theta/4\omega)t}$and taking above initial condition the expression of jacobian takes following form 
\begin{align}
\xi_{\eta}=\frac{(1-e^{-\theta(2+\gamma/4)t})3cos(\eta)}{\theta(8+\gamma)}
\left[4\alpha\left(\frac{2}{\gamma \theta}(1-\beta)+sin(\eta)\right)-
(1+\gamma)\beta\left(\frac{\theta}{(\gamma-4)}-\frac{2(2-\gamma)}{\gamma\theta}\right)\right]\nonumber\\
+\frac{(1-e^{-\theta(1+\gamma/2)t})\beta cos(\eta)}{2\theta(2+\gamma)}
\left[-3(2+\gamma)\left(\frac{2}
{\gamma \theta}+sin(\eta)\right)+\frac{\gamma}{4}\right]\qquad\qquad\qquad\qquad\qquad\qquad\nonumber\\
+\frac{(1-e^{-\theta\gamma t/4})(1+\gamma) cos(\eta)}{\gamma\theta}
\left[\frac{3(4-\gamma\beta)}{\gamma \theta }+\frac{3\gamma \beta}{\theta(\gamma-4)}+
2\left(1-\frac{\gamma \beta}{4}\right)\right]\qquad\qquad\qquad\qquad\nonumber\\
-\frac{(1-e^{-\theta\gamma t/2})3(1+\gamma) cos(\eta)}{2\gamma\theta}
\left[\frac{(4-\gamma \beta)}{\gamma \theta}+2sin(\eta)\left(1-\frac{\gamma \beta}{2}+
\right)\right]+1\qquad\qquad\qquad\qquad\nonumber\\
\end{align}

\subsection{Bounds on parameters}
In our case we have  choosen $ \Gamma $ such that
\begin{equation}\label{l1.46}
\Gamma = \frac{(\gamma+1)}{2}\left(\frac{1}{(1-\beta\rho_{0})}-
\frac{2\alpha\rho_{0}(2-\gamma-3\beta\rho_{0})}{a_{0}^{2}(\gamma+1)(1-\beta\rho_{0})}\right)|\nabla\phi|
=O(\epsilon),
\end{equation}
also substituting the corresponding expression for $\rho_{0}  $ and $ a_{0} $ the expression of $ \Gamma  $ can be  rewritten into the form
\begin{equation}\label{l1.47}
\frac{(\gamma+1)}{2}\left(1-
\frac{2\alpha(1+\omega t)^{2-\theta/\omega}(2-\gamma-3\beta(1+\omega t)^{-\theta/\omega})}{a_{0}^{2}(\gamma+1)}\right)\frac{|\nabla\phi|}{(1-\beta(1+\omega t)^{-\theta/\omega})}.
\end{equation}
In order to choose $ \Gamma=O(\epsilon) $ along with the conditions $ \alpha ,$  $ \beta ,$ $ \omega ,$ $ \theta $  all are positive also $ t \ge 0 $ and $ 1< \gamma\leq 5/3 $. If we assume $ \theta = \omega ,$ then We obtained the following expression  relating $ \alpha ,$ $ \beta, $ and $ t $
\begin{equation}\label{1.48}
t=\dfrac{1}{\omega} \left[\left(\frac{(\gamma +1)}{2\alpha}+3\beta\right)\dfrac{1}{(2-\gamma)}-1\right],
\end{equation}
which in view of $ t \geq 0 $ gives the relation $ \dfrac{(\gamma +1)}{2 \alpha} + 3\beta \geq (2-\gamma).$ We have taken $ \gamma=1.01 $ and $ \omega=0.1 $ and hence the above condition reduces to 
\begin{equation}\label{l1.48}
\left(\frac{67}{66-200\beta}\right)\geq \alpha \qquad\text{and } \qquad
t_{0}=10\left[\left(\frac{67}{200\alpha}+\beta\right)\dfrac{100}{33}-1\right],
\end{equation}
we have used these relation in our numerical calculation for finding the values of $ \alpha $ and $ \beta, $ and the initial time $ t_{0} $ for corresponding $ \Gamma = O(\epsilon) .$

\section{Effects of parameters on breaking of solution}\label{real1}
In this section we have shown the effect of parameters $ \alpha, $ $ \beta, $ $ \theta, $ $ \omega $ on breaking of solution. We have taken initial data

\begin{equation}\label{examplesin}
\sigma(\xi,t_{0})=
\left\{ \begin{array}{l l}
sin(\xi) & \quad\ \xi\in (0,\pi ),\\
~~0 & \quad \text{~otherwise}
\end{array} \right.
\end{equation} 
where $ t_{0} $ is the value obtain of time obtained from the restriction $ \Gamma =O(\epsilon). $

\subsubsection{Effects of $\beta$ on breakdown of solution:} 
\begin{figure}[!t]
	\centering
	\includegraphics[width=0.9\textwidth]{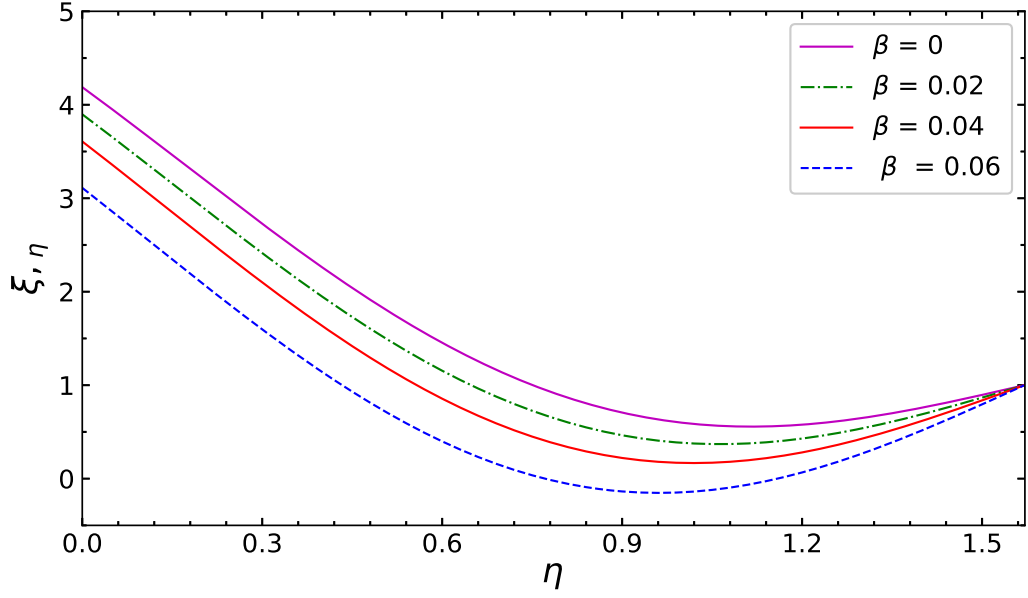}
	\caption{Graph of the jacobian $\xi_{,\eta}$ for various values of $\beta$ with $\gamma=1.01$,
		$\delta=0.1$, $\omega=0.1$, $ \alpha=0.35 $ at $ t=1.4 $.}
	\label{fluxb07}
\end{figure}

To discuss the effect of $ \beta $ on the breaking of solution we have varied the value of $ \beta $ while keeping all the other parameters as constant with values $ \alpha=0.15 ,$ $ \gamma=1.01, $ $ \theta=0.1, $ $\omega=0.1 , $ $ t=1.4 ,$ and plotted the Jacobian $ \xi_{,\eta} $ against $ \eta $ and found that with the increase in beta nonlinear effect serve to expedite as seen in the Figure \ref{fluxb07} .
\subsubsection{Effects of $ \alpha $ on breakdown of solution:} 
\begin{figure}[!t]
	\centering
	\includegraphics[width=0.9\textwidth]{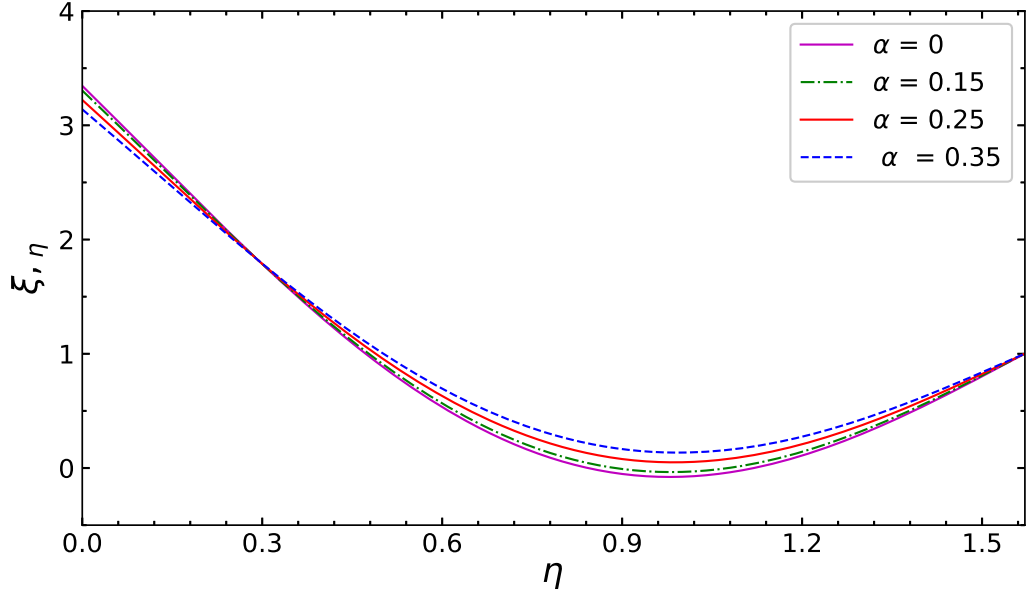}
	\caption{Graph of the jacobian $\xi_{,\eta}$ for various values of $\alpha$ with $\gamma=1.01$,
		$\delta=0.1$, $\omega=0.1$, $ \beta=0.06 $ at $ t=1.4 $.}
	\label{fluxb007}
\end{figure}

All the parameter with same values are taken as in the last case except $ \alpha $ and the effects of the variation in  $ \alpha $ are noticed keeping all other parameters as constants. In contrast to the last case, with the increase in $ \alpha  $ nonlinear effect serves to delay as shown in Figure \ref{fluxb007} .

\subsubsection{Effects of $ \theta $ and $ \omega$ on breakdown of solution:} 
\begin{figure}[!t]
	\centering
	\includegraphics[width=0.9\textwidth]{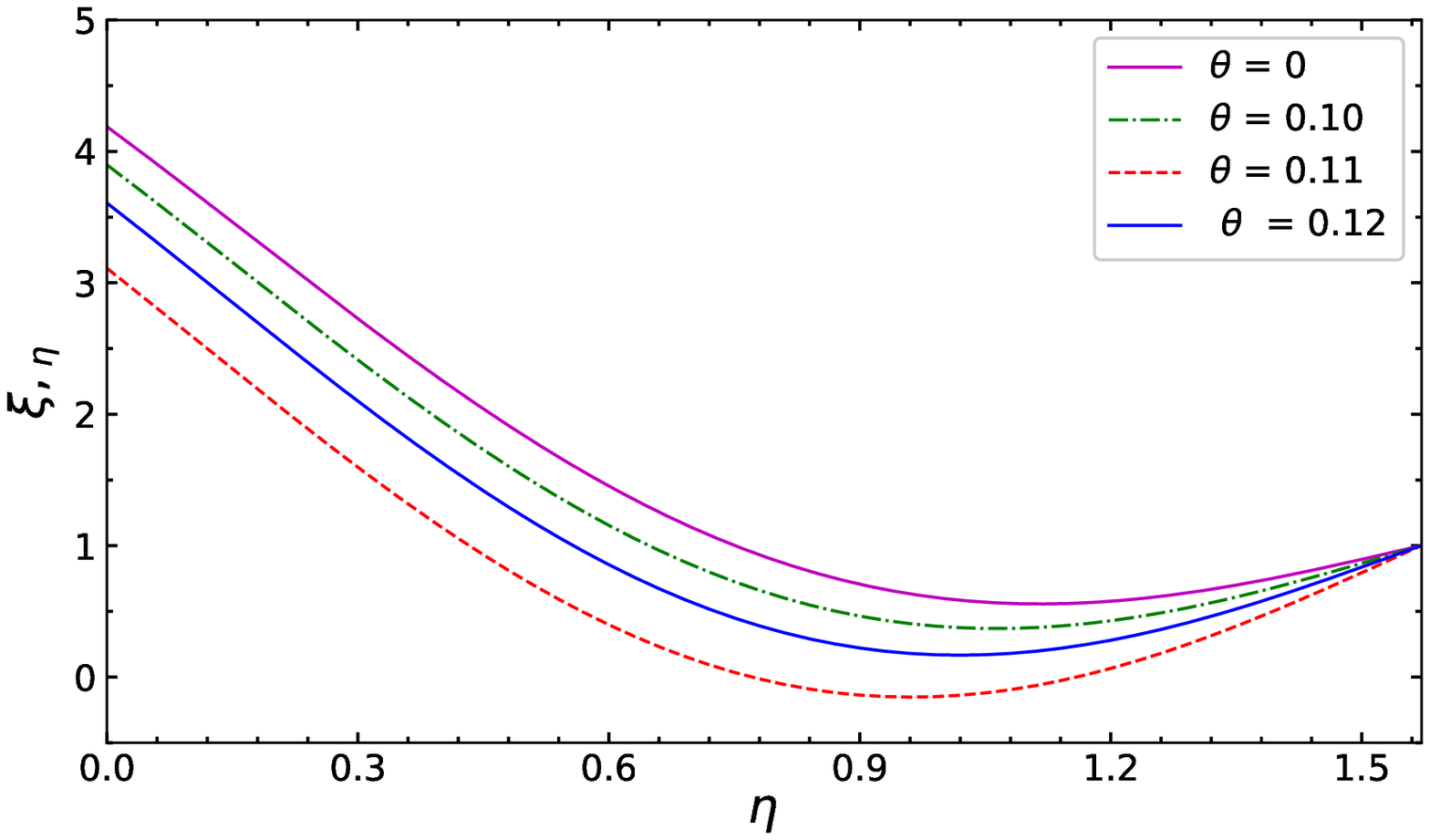}
	\caption{Graph of the jacobian $\xi_{,\eta}$ for various values of $\theta$ with $\gamma=1.01$,
		$\delta=0.1$, $\omega=0.1$, $\alpha=0.35,$ $ \beta=0.06 $ at $ t=1.4 $.}
	\label{fluxb0007}
\end{figure}

\begin{figure}[!t]
	\centering
	\includegraphics[width=0.9\textwidth]{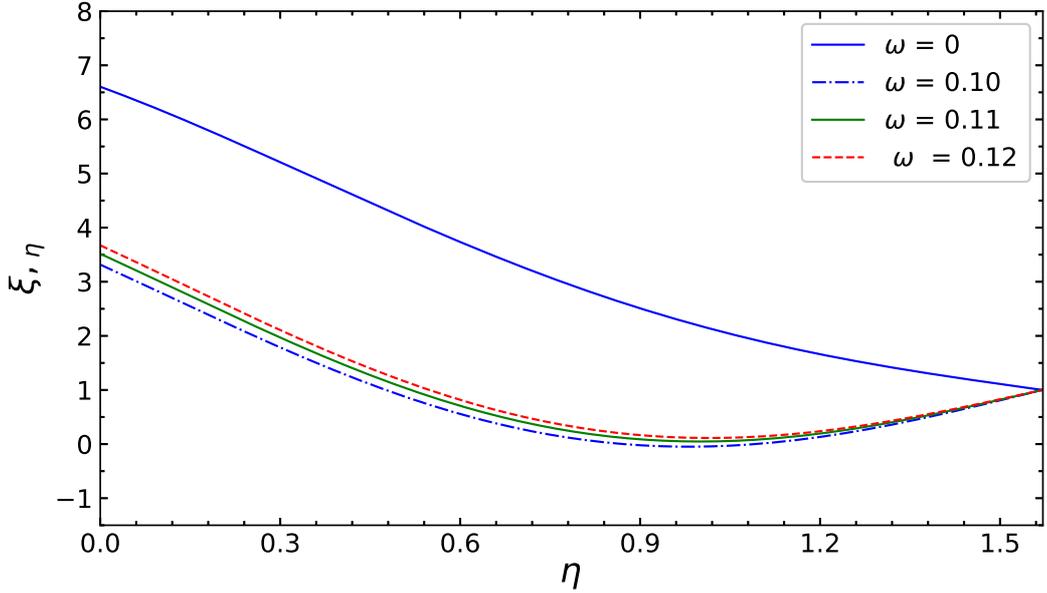}
	\caption{Graph of the jacobian $\xi_{,\eta}$ for various values of $\omega$ with $\gamma=1.01$,
		$\delta=0.1$, $\alpha=0.35$, $ \beta=0.06 $ at $ t=1.4 $.}
	\label{fluxb00007}
\end{figure}

Finally, the effects of the atmospheric parameters, i.e., density variation parameter $ \theta  $ and sound speed variation parameter $\omega $ are observed. It is find out that, increase in the density parameter helps in the breakdown whereas speed of sound variation parameter delays breaking of solution as displayed in Figs. \ref{fluxb0007},  \ref{fluxb00007}, respectively.

\subsection{Evolution of waves in a van der Waals fluid}
\begin{figure}[!t]
	\centering
	\includegraphics[width=1.00\textwidth]{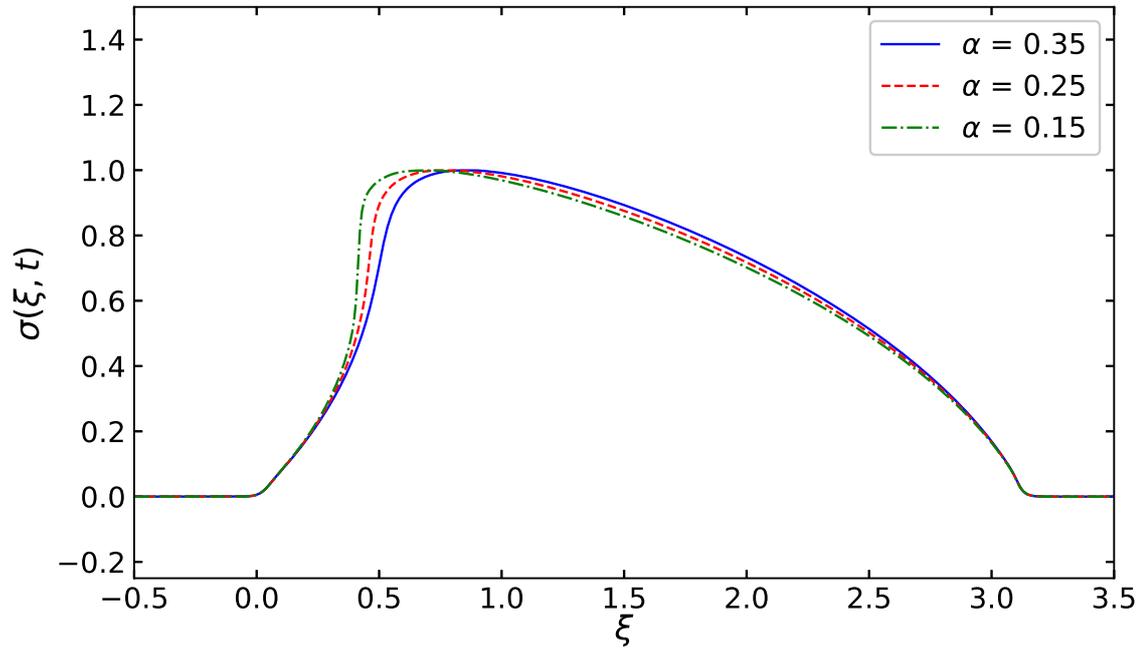}
	\caption{Numerical solutions of \eqref{l1.42} and \eqref{examplesin} with $\gamma=1.01$, $\beta=0.06$, and $\theta=0.1$, $\omega=0.1$, at time $t=1.4$; for $\alpha=0.15$, $\alpha=0.25$, and  $\alpha=0.35$.}
	\label{KineEx1}
\end{figure}

\begin{figure}[!t]
	\centering
	\includegraphics[width=1.00\textwidth]{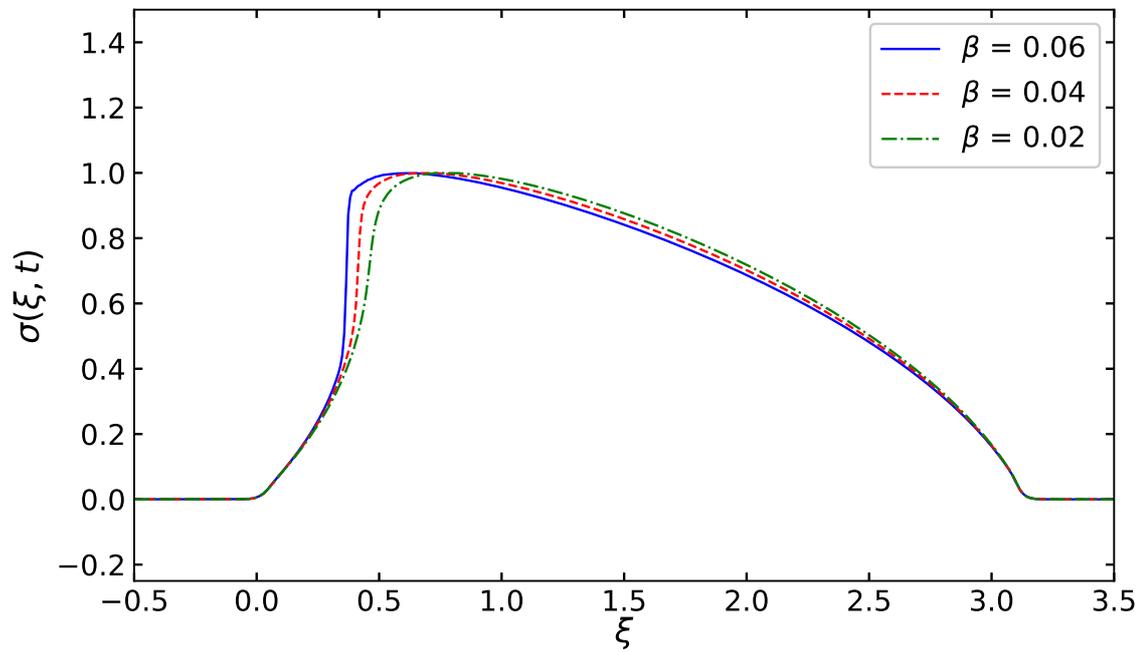}
	\caption{Numerical solutions of \eqref{l1.42} and \eqref{examplesin} with  $\gamma=1.01$, $\alpha=0.35$, $\theta=0.1$, and  $\omega=0.1$, at time $t=1.4$; for $\beta=0.02$,  $\beta=0.04$, and  $\beta=0.06$.}
	\label{KineEx2}
\end{figure}

In this section, to discuss the effects of van der Waals parameters ($ \alpha, \,\beta $) in the case when flux function of evolution equation has  quadratic as well cubic nonlinearity, we present numerical solution of evolution equation 
$ (\ref{l1.42}) $ with the following initial data:
\begin{equation}\label{examplesin1}
\sigma(\xi,t_{0})=
\left\{ \begin{array}{l l}
sin(\xi) & \quad\ \xi\in (0,\pi),\\
~~0 & \quad \text{~otherwise}.
\end{array} \right.
\end{equation} 
For various values of $ (\alpha,\, \beta )$,  other parameters are so chosen such that the condition (\ref{l1.46}), i.e. 
$ \Gamma = O(\epsilon) $ is satisfied, we observe that the breaking of solution delays  with an increase in  $ \alpha $ as is exhibited in Figure
\ref{KineEx1}; however, an increase in $ \beta $ has  just the opposite effect, i.e., the breaking
of solutions gets expedited  as seen in the Figure  \ref{KineEx2}. Here 
we  notice that it is  the cubic nonlinearity in the flux function that is responsible for the  breakdown of solution on an expansion phase of the wave profile, which is quite different from the quadratic flux case, where there is no breakdown on the expansion phase if the initial datum is monotonic increasing. 

%

\newpage
\section{Conclusions}\label{conv}
We have studied, using perturbation methods, propagation of high frequency waves with mixed nonlinearity in a stratified atmosphere with van der Waals equation of state. A transport equation for the wave amplitude is derived; which exhibits  both quadrartic and cubic nonlinearities. A quiet and steady atmosphere with thermodynamic quantities, depending  only on one spatial coordinate (height) with varying density and sound speed, is considered. It is shown that the Cauchy problem   exhibits a breakdown of the continuous solution on the expansive phase of the wave profile, which is monotonic increasing, in the sense that the Jacobian of the transformation vanishes after a finite time.  This behaviour is due to the presence of the cubic nonlinearity term in the flux function and is quite different from   the quadratic nonlinearity case where the solution is always continuous.
Effects of the influence of  van der Waals parameters $ \alpha,\, \beta  $ on the breaking of solution is displayed in Figs. \ref{fluxb07}, \ref{fluxb007}, respectively, while effects of atmospheric parameters $ \theta,\, \omega  $ was observed in Figs. \ref{fluxb007} and \ref{fluxb0007}. Indeed, the effect of the van der Waals parameter $ \alpha $ is to delay the onset of singularity in the solution, whereas the effect of $ \beta $ is to hasten the process of singularity formation in the solution as shown in Figs. 
\ref{KineEx1}, \ref{KineEx2}.



\section*{References}

\begin{thebibliography}{99}
	
	\bibitem{MR1166185}	M. S. Cramer and R. Sen. A general scheme for the derivation of evolution equations describing mixed nonlinearity. Wave Motion, 15(4):333-355, 1992.
	
	\bibitem{MR716196} J. K. Hunter and J. B. Keller. Weakly nonlinear high frequency waves. Comm. Pure
	Appl. Math., 36(5):547-569, 1983.
	
    \bibitem{MR1600236} A. Kluwick and E. A. Cox. Nonlinear waves in materials with mixed nonlinearity.
    Wave Motion, 27(1):23-41, 1998.
    
	\bibitem{Lax} P. D. Lax. Hyperbolic Systems of Conservation Laws and the Mathematical Theory of
	shock Waves. SIAM, Philadelphia, 1973
	
	\bibitem{MR760229} A. Majda and R. Rosales. Resonantly interacting weakly nonlinear hyperbolic waves. I. A single space variable. Stud. Appl. Math., 71(2):149-179, 1984.
	 
	\bibitem{NZ73} J. W. Nunziato and E. K. Walsh. Shock-wave propagation in inhomogeneous atmo-
	spheres. Physics of Fluids, 16:482-484, 1973.
			
	\bibitem{JF85} J. A. Sanders and F. Verhulst. Averaging methods in nonlinear dynamical systems.
	Springe-Verlag, New York, 1985.

   \bibitem{zhao2011admissible} N. Zhao, A. Mentrelli, T. Ruggeri, and M. Sugiyama. Admissible shock waves and shock-induced phase transitions in a van der waals 
uid. Physics of 
fluids,
   23(8):086101, 2011.
	
	
    
\end{thebibliography}

\end{document}